\documentclass[11pt]{amsart}

\usepackage{euscript,oldgerm,geometry}
\usepackage{amssymb,amsmath}
\usepackage{color} 
 
\newcommand{\conju}[1]{\overline{#1}}
\newcommand{\refe}[1]{(\ref{#1})}
\newcommand{\dst}{\displaystyle}

\newcommand{\NN}{{\mathbb N}}

\newcommand{\CC}{{\mathbb C}}

\renewcommand{\u}{\mbox{${\textswab u}$}}
\newcommand{\newu}{\mbox{$\widetilde{\textswab u}$}}
\newcommand{\p}{\widetilde{P}}

\newcommand{\K}{\mathrm{K}}
\newcommand{\Ku}{\mathrm{K}^{\alpha,\beta}}

\newcommand{\bq}{\begin{equation}}
\newcommand{\eq}{\end{equation}}
\newcommand{\ba}{\begin{array}}
\newcommand{\ea}{\end{array}}

\newtheorem{theorem}{Theorem}
\newtheorem{rem}[theorem]{Remark}

\begin{document}

\title{The $q$-Racah-Krall-type polynomials}

\author{R. \'{A}lvarez-Nodarse}
\address{IMUS \& Departamento de An\'alisis Matem\'atico, Universidad de Sevilla.
Apdo. 1160, E-41080  Sevilla, Spain}\email{ran@us.es}

\author{R. Sevinik Ad\i g{\"{u}}zel}
\address{Sel\c{c}uk University, Faculty of Science, Department of Mathematics, 42075,
Konya, Turkey}\email{sevinikrezan@gmail.com}

\thanks{The first author was partially supported by DGES grants MTM2009-12740-C03;
PAI grant FQM-0262 and Junta de Andaluc\'\i a grant P09-FQM-4643.}



\begin{abstract}
In this paper the Krall-type polynomials obtained 
via the addition of two mass points to the weight function of the 
\textit{standard} $q$-Racah polynomials are introduced. Several algebraic 
properties of these 
polynomials are obtained and some of their limit cases are discussed.
\end{abstract}

\maketitle

\section{Introduction}

In the last decades the study of the discrete analogues of the classical 
special functions and, in particular, the orthogonal polynomials,
has received an increasing interest (for a review see \cite{ran,ks,NSU}).
Special emphasis was given to the {\em q-}analogues of the orthogonal polynomials or {\em q-}polynomials, 
which are closely related with different topics in other fields of actual science: 
Mathematics and Physics (see e.g., \cite{aac1,aac2,NSU,koor-tut,vk} and references therein).
One of the possible extensions of the classical polynomials are the so-called 
Krall-type polynomials.

The Krall-type polynomials are polynomials which are orthogonal with respect to a
linear functional $\newu$ obtained from a
quasi-definite functional $\u:\mathbb{P}\mapsto\CC$  ($\mathbb{P}$, denotes the space
of complex polynomials with complex coefficients) via the addition of delta Dirac measures.
In the last years the study of such polynomials have attracted an increasing interest
with a special emphasis on the case when the starting functional $\u$ is a classical
\textit{continuous}, discrete or $q$ linear functional (for more details see  \cite{RR} 
and references therein). These kind of polynomials appear as eigenfunctions of a fourth order 
linear differential operator with polynomial coefficients that do not depend on the degree of the polynomials. 
They were firstly considered by Krall in 1940 (see e.g. \cite[chapter XV]{akrall}) and
further studied {by} several authors (for more details see  \cite{RR} and references therein).
For the case of the discrete lattice A. Dur\'an {has} discovered very recently 
\cite{dur1,dur2} a method for obtaining the orthogonal polynomials satisfying 
higher order differential and difference equations. 

In two recent papers \cite{RR,RR1} a general theory of 
{the} Krall-type polynomials on non-uniform lattices was developed. In fact, in 
\cite{RR,RR1} the authors studied {the} polynomials $\p_n(s)_q$ which are orthogonal with respect
to the linear functionals $\newu=\u+\sum_{k=1}^N A_k\delta_{x_k}$ defined on the $q$-quadratic lattice 
$x(s)=c_1q^s+c_2q^{-s}+c_3$. For these polynomials the following expression have
been found \cite[Eq. (13)]{RR}
\begin{equation}
\p_n(s)_q=P_n(s)_q-\sum_{i=1}^M A_i\p_n(a_i)_q\K_{n-1}(s,a_i),
\label{repform}
\end{equation}
provided that the sequence of monic polynomials $(\p_n)_n$ exist.
Here, $(P_n)_n$ are the polynomials orthogonal with respect to $\u$,
$A_i\in\mathbb{C}$ and $K_n(x,y)$ are the $n$-th
reproducing kernels given by the Christoffel-Darboux formula
\begin{equation*}\label{kern}
\begin{array}{l}
\K_n(s_1, s_2):=\dst\sum_{k=0}^{n}{\dst\frac{P_k(s_1)_qP_k(s_2)_q}{d_k^2}}=\dst\frac{\alpha_n}{d_n^2}
\frac{P_{n+1}(s_1)_q P_n(s_2)_q-P_{n+1}(s_2)_qP_{n}(s_1)_q}
{x(s_1)-x(s_2)}.
\end{array}
\end{equation*}
Moreover, in \cite{RR}, we have studied 
the modifications of the \textit{non-standard} $q$-Racah polynomials
defined on the lattice $x(s)=[s]_q[s+1]_q$ \cite{RYR,NSU}, where $[s]_q$ denotes the 
symmetric $q$-numbers $[s]_q=({q^{s/2}-q^{-s/2}})/({q^{1/2}-q^{-1/2}})$,
whereas in \cite{RR1} we developed the theory for the Askey-Wilson-Krall polynomials.
Here, in this paper we will consider the ``\textit{standard}'' (or classical)
$q$-Racah-Krall polynomials defined on the lattice 
$x(s)=q^{-s}+\delta q^{-N}q^s$ introduced in \cite{AW} (see also 
\cite[page 422]{ks}). For doing that we will follow the 
general results obtained in \cite[section 2, 3]{RR}.

The structure of the paper is as follows. In Section 2, 
the \textit{standard} $q$-Racah-Krall polynomials are introduced and 
all their main characteristics are studied in detail.
Finally in  Section 3, some important limit cases are considered.

\section{Standard $q$-Racah-Krall polynomials}

In this section, we study the \textit{standard} $q$-Racah-Krall polynomials, i.e., the
polynomials obtained from the \textit{standard} $q$-Racah {polynomials} $R_n^{\alpha,\beta}(s)_q$ via the addition
of two mass points. In the following we follow the notations introduced in \cite{RR}.

\subsection{Preliminaries}
The \textit{standard} $q$-Racah $R_n^{\alpha,\beta}(s)_q:= R_n(x(s),\alpha,\beta, \gamma, \delta|q)$
are polynomials on the $q$-quadratic lattice $x(s)=q^{-s}+\delta \gamma q^{s+1}$ introduced 
in  \cite{AW} (see also the  more recent book \cite[page 422]{ks}). They are 
defined by the following basic series (for properties of basic series see \cite{GR}) 
\bq
\begin{split}\dst
R_n(x(s),\alpha,\beta, \gamma, \delta|q) & = \dst
\frac{(\alpha q, \beta\delta q, \gamma q; q)_n }{(\alpha\beta q^{n+1};q)_n} 
{}_{4}\varphi_3 \left(\ba{c} q^{-n}, \alpha\beta q^{n+1}, q^{-s},
\delta \gamma q^{s+1}  \\ \alpha q, \beta\delta q, \gamma q  \ea
\bigg| q ; q \right),
\end{split}
\label{pol-rac}
\eq 
where $x(s)=q^{-s}+\delta \gamma q^{s+1}$,
and $\alpha q=q^{-N}$, $\beta\delta q=q^{-N}$,  $\gamma q =q^{-N}$, $N$, a nonnegative integer.
Their main data can be found in table \ref{tabla1}.
Notice from \refe{pol-rac} that we have
\begin{equation*}
\begin{split}
R_n^{\alpha,\beta}(0)_q:=R_n(x(0),\alpha,\beta, \gamma, \delta|q) & = \dst\frac{(\alpha q, \beta\delta q, \gamma q; q)_n}
{(\alpha\beta q^{n+1}; q)_n},
\end{split}
\end{equation*}
and $R_n^{\alpha,\beta}(N)_q:=R_n(x(N),\alpha,\beta, \gamma, \delta|q)$
\begin{equation*}
\begin{split}
R_n^{\alpha,\beta}(N)_q & =(\delta\gamma q^{N+1})^n\frac{(q^{-N}, \beta\gamma^{-1} q^{-N}, \delta^{-1} q^{-N}; q)_n}
{(\beta q^{n-N}; q)_n},\quad \mbox{if } \alpha q=q^{-N},\\
R_n^{\alpha,\beta}(N)_q&  =\delta^n\frac{(q^{-N}, \beta q, \delta^{-1}\alpha q; q)_n}
{(\alpha\beta q^{n+1}; q)_n}, \quad \mbox{if } \gamma q=q^{-N},\\
R_n^{\alpha,\beta}(N)_q& =(\beta^{-1}\gamma)^n\frac{(q^{-N}, \beta q, \delta^{-1}\gamma^{-1}\alpha q^{-N}; q)_n}
{(\alpha\beta q^{n+1}; q)_n}, \quad \mbox{if } \beta\delta q=q^{-N}.
\end{split}
\end{equation*}

\begin{table}[ht!]
\caption{Main data of the monic $q$-Racah polynomials 
\label{tabla1}}
\begin{center}
\begin{scriptsize}
{\renewcommand{\arraystretch}{.25}
\begin{tabular}{|@{}c@{}| | @{}c@{}|}\hline
 &  \\
$P_n(s)$ & $R_n(x(s),\alpha,\beta,\gamma,\delta|q)   \,,\quad x(s) = q^{-s}+\delta \gamma q^{s+1}$ , 
\,\,\, $\Delta x(s)=q^{-s}(1-\delta\gamma q^{2s+2})(q^{-1}-1)$ \\
 &  \\
\hline\hline
 &  \\
$ (a,b)$ & \mbox{
$[0,N]$} \\
 &  \\
\hline
 &  \\
$\rho(s)$ &  \mbox{$\dst\frac{(\alpha\beta)^{-s}}{(q^{-1/2}-q^{1/2})}
\frac{(\alpha^{-1}\beta^{-1}q^{-1}, \alpha^{-1}\gamma\delta q,\beta^{-1}\gamma q,\delta q; q)_{\infty}}
{(\alpha^{-1}\beta^{-1}\gamma, \alpha^{-1}\delta, \beta^{-1}, \gamma\delta q^2; q)_{\infty}}
\frac{(\delta \gamma q, \alpha q, \beta\delta q, \gamma q; q)_s}
{(q, \alpha^{-1}\delta \gamma q, \beta^{-1}\gamma q, \delta q; q)_s}$}  \\
 &  \\
\hline
 &  \\
$\sigma(s)$ & \mbox{
$  \delta^2q^{-2N}(q^{1/2}-q^{-1/2})^2q^{-2s}(q^s-1)(q^s-\delta^{-1})
(q^s-\beta \gamma^{-1})(q^s-\alpha\delta^{-1}\gamma^{-1}) $} \\
 &  \\
\hline
 &  \\
$\Phi(s)$ & \mbox{
$  (q^{1/2}-q^{-1/2})^2q^{-2s}(1-\alpha q^{s+1})(1-\beta\delta q^{s+1})(1-\gamma q^{s+1})(1-\delta \gamma q^{s+1}) $} \\
 &  \\
\hline
 &  \\
$\tau(s) $ &  \mbox{
\begin{tabular}{r}
$ \frac{(q^{-1/2}-q^{1/2})q^{-s}}{(1-\gamma\delta q^{2s+1})}\Big[(1-\alpha q^{s+1})(1-\beta\delta q^{s+1})(1-\gamma q^{s+1})(1-\delta \gamma q^{s+1})-(\delta\gamma q)^2 (q^s-1)$  \\[3mm]
 $\times (q^s-\delta^{-1})(q^s-\beta \gamma^{-1})(q^s-\alpha\delta^{-1}\gamma^{-1})\Big]$ \end{tabular}} \\
 &  \\
 \hline
$\tau_n(s) $ &  \mbox{
\begin{tabular}{r}
$ -q^{-n}(q^{1/2}-q^{-1/2})\Big\{(1-\alpha\beta q^{2n+2})x(s+\frac n2)+q^{-n/2}\Big[(1-\alpha q^{n+1})
(1-\beta\delta q^{n+1})(1-\gamma q^{n+1})
$\\[3mm] $-(1+\delta \gamma q^{n+1})(1-\alpha\beta q^{2n+2})\Big]\Big\}$ \end{tabular}}\\
 &  \\
 \hline
 &  \\
$\lambda _n$ &\mbox{
$-q^{-n+\textstyle \frac 12}(1-q^n)(1-\alpha\beta q^{n+1})$} \\
 &  \\
\hline
 &  \\
$d_n^2$ &    \mbox{$\dst\frac{(1-\alpha\beta q)(\delta \gamma q)^n}{(1-\alpha\beta q^{2n+1})}\dst\frac{(\alpha q, \beta\delta q, \gamma q, q, q\alpha\beta \gamma^{-1},
\alpha\delta^{-1}q, \beta q; q)_n}{(\alpha\beta q,\alpha\beta q^{n+1},\alpha\beta q^{n+1}; q)_n}
$}\,\,\, \\
 &  \\
 \hline
 $\beta_n $ &  \begin{tabular}{c}
$\dst1+\delta \gamma q-\frac{(1-\alpha q^{n+1})(1-\alpha\beta q^{n+1})(1-\beta\delta q^{n+1})(1-\gamma q^{n+1})}
{(1-\alpha\beta q^{2n+1})(1-\alpha\beta q^{2n+2})}$ \\[3mm]
$\dst-\frac{q(1-q^n)(1-\beta q^n)(\gamma-\alpha\beta q^n)(\delta-\alpha q^n)}
{(1-\alpha\beta q^{2n})(1-\alpha\beta q^{2n+1})}$
\end{tabular}\\
 &  \\
\hline
 &  \\
$\gamma_n$ & $ \dst\frac{(1-\alpha q^{n})(1-\alpha\beta q^{n})(1-\beta\delta q^{n})(1-\gamma q^{n})}
{(1-\alpha\beta q^{2n-1})(1-\alpha\beta q^{2n})}$ 
$\dst\frac{q(1-q^n)(1-\beta q^n)(\gamma-\alpha\beta q^n)(\delta-\alpha q^n)}
{(1-\alpha\beta q^{2n})(1-\alpha\beta q^{2n+1})} $ \\
 &  \\
\hline
&  \\
$\overline{\alpha}_n=\widehat{\alpha}_n$ &  $q^{-n+\textstyle\frac 12}(q^{-1/2}-q^{1/2})(1-\alpha\beta q^{2n+1})$  \\&  \\
\hline &  \\
$\overline{\beta}_n(s) $ & $\begin{tabular}{r}
$q^{-n/2+1/2}(q^{1/2}-q^{-1/2})\dst\frac{(1-\alpha\beta q^{n+1})}{(1-\alpha\beta q^{2n+2})}\Big\{(1-\alpha\beta q^{2n+2})x(s+\frac n2)+q^{-n/2}\Big[(1-\alpha q^{n+1})
$\\[3mm] $\times (1-\beta\delta q^{n+1})(1-\gamma q^{n+1})-(1+\delta \gamma q^{n+1})(1-\alpha\beta q^{2n+2})\Big]\Big\}$ \end{tabular}$\\
\hline &\\
$\widehat{\beta}_n(s) $ & 
$\overline{\beta}_n(s)-q^{-s-n+\textstyle\frac 12}(q^{1/2}-q^{-1/2})(1-q^n)(1-\alpha\beta q^{n+1})(1-\delta \gamma q^{2s+1})$ \\
\hline
\end{tabular}  }
\end{scriptsize}
\end{center}
\end{table}
We also need the following identity for the \textit{standard} $q$-Racah polynomials
which follows from \cite[Eq. (7)]{RR}
\bq\label{ThetaXi}
R_{n-1}^{\alpha,\beta}(s)_q=\Theta(s,n)R_n^{\alpha,\beta}(s)_q+\Xi(s,n) R_n^{\alpha,\beta}(s+1)_q,
\eq
where
\begin{footnotesize}
\begin{equation*}
\begin{split}
&\Theta(s,n)=-\Xi(s,n)+
\frac{(1-\alpha\beta q^{2n-1})(1-\alpha\beta q^{2n})^2
}{(1-\alpha q^n)(1-\alpha\beta q^{n})(1-\beta\delta q^n)(1-\gamma q^{n})(1-q^n)(1-\beta q^n)}\\
&\times\frac{1}{(\gamma-\alpha\beta q^n)(\delta-\alpha q^n)}
\Big[(q^{-s}-1)(1-\delta \gamma q^{s+1})
+\frac{(1-\alpha q^{n+1})(1-\alpha\beta q^{n+1})}
{(1-\alpha\beta q^{2n+1})(1-\alpha\beta q^{2n+2})}\\
&\times(1-\beta\delta q^{n+1})(1-\gamma q^{n+1})\dst+\frac{q(1-q^n)(1-\beta q^n)(\gamma-\alpha\beta q^n)(\delta-\alpha q^n)}
{(1-\alpha\beta q^{2n})(1-\alpha\beta q^{2n+1})}\Big]\\
&-\frac{(1-\alpha\beta q^{2n-1})(1-\alpha\beta q^{2n})^2
}{(1-\alpha q^n)(1-\alpha\beta q^{n})(1-\beta\delta q^n)
(1-\gamma q^{n})(1-q^n)(1-\beta q^n)(\gamma-\alpha\beta q^n)}\\
&\times\frac{1}{(\delta-\alpha q^n)}\Bigg\{q^{n/2}\dst
\frac{(1-\alpha\beta q^{n+1})}{(1-\alpha\beta q^{2n+2})}\Big\{(1-\alpha\beta q^{2n+2})x_n(s)+q^{-n/2}\Big[(1-\alpha q^{n+1})\\
&\times(1-\beta\delta q^{n+1})(1-\gamma q^{n+1})-(1+\delta\gamma q^{n+1})(1-\alpha\beta q^{2n+2})\Big]\Big\}-q^{-s}(1-q^n)\\
&\times(1-\alpha\beta q^{n+1})(1-\delta \gamma q^{2s+1})\Bigg\},
\end{split}
\end{equation*}
\end{footnotesize}

\begin{footnotesize}
\begin{equation*}
\begin{split}
&\Xi(s,n)=\dst-\frac{q^{-s-1+n}(1-\alpha\beta q^{2n-1})(1-\alpha\beta q^{2n})^2
(1-\alpha q^{s+1})}
{(1-\alpha q^n)(1-\alpha\beta q^{n})(1-\beta\delta q^n)(1-\gamma q^{n})(1-q^n)(1-\beta q^n)}\\
&\times\frac{(1-\beta\delta q^{s+1})(1-q^{-N+s})(1-\delta\gamma q^{s+1})}
{(\gamma-\alpha\beta q^n)(\delta-\alpha q^n)(1-\delta \gamma q^{2s+2})}.
\end{split}
\end{equation*}
\end{footnotesize}


\subsection{The case of two mass point}\label{2masspoint}

Let us now consider the modification of the \textit{standard} $q$-Racah polynomials defined
in \refe{pol-rac}. Let $\u$ be the functional with respect to which the \textit{standard} $q$-Racah 
polynomials \refe{pol-rac} are orthogonal. Then we study the orthogonal polynomials
that are orthogonal with respect to the following linear functional 
$\newu=\u+A\delta(x(s)-x(0))+B\delta(x(s)-x(N))$. In other words, the $q$-Racah-Krall 
polynomials, ${R}_n^{\alpha,\beta,A,B}(s)_q:={R}_n^{A,B}(x(s), \alpha, \beta, \gamma, \delta|q)$ 
satisfying the following orthogonality relation
\bq\label{ort-umod}
\begin{split}
\sum_{s = 0 }^{N} {R}_n^{\alpha,\beta,A,B}(s)_q {R}_m^{\alpha,\beta,A,B}(s)_q 
\rho(s)&\Delta x(s-1/2) + A {R}_n^{\alpha,\beta,A,B}(0)_q
{R}_m^{\alpha,\beta,A,B}(0)_q \\ & +B{R}_n^{\alpha,\beta,A,B}(N)_q
{R}_m^{\alpha,\beta,A,B}(N)_q = \delta_{n,m} \widetilde{d}_n^2,
\end{split}
\eq
where $\rho$ is the \textit{standard} $q$-Racah weight function (see table
\ref{tabla1} (as in \cite{RR} we have chosen $\rho(s)$ to be a probability measure).
From \cite[Eqs. (31)-(32)]{RR} it follows that 
\begin{equation}\label{u-u2m}
\begin{array}{l}\dst
{R}_n^{\alpha,\beta,A, B}(0)_q=\frac{(1+B \Ku_{n-1}(N, N))R_n^{\alpha,\beta}(0)_q-
B\Ku_{n-1}(0, N)R_n^{\alpha,\beta}(N)_q}{\kappa_{n-1}^{\alpha,\beta}(0,N)},\\
{R}_n^{\alpha,\beta,A, B}(N)_q=\dst
\frac{-A\Ku_{n-1}(N, 0)R_n^{\alpha,\beta}(0)_q+(1+A\Ku_{n-1}(0, 0))R_n^{\alpha,\beta}(N)_q}
{\kappa_{n-1}^{\alpha,\beta}(0,N)},
\end{array}
\end{equation}\begin{small}
\begin{equation*}
\begin{split}
\widetilde{d}_n^2  =
d_n^2 &\! + \!\frac{A ({R}_n^{\alpha,\beta}(0)_q)^2 \{1+ B\Ku_{n-1}(N,N)\}  +
B({R}_n^{\alpha,\beta}(N)_q)^2 \{ 1+ A\Ku_{n-1}(0,0) \} }{\kappa_{n-1}^{\alpha,\beta}(0, N)} \\
& -\frac{2 AB {R}_n^{\alpha,\beta}(0)_q {R}_n^{\alpha,\beta}(N)_q\Ku_{n-1}(0,N)}
{\kappa_{n-1}^{\alpha,\beta}(0, N)},
\end{split}
\end{equation*}\end{small}
where
\begin{equation}\label{K(a,b)}
\begin{split}
\kappa_{m}^{\alpha,\beta}(s,t)& =1+A\Ku_{m}(s,s)+ B\Ku_{m}(t,t)\\
& \quad +AB\left\{\Ku_{m}(s,s)\Ku_{m}(t,t)-(\Ku_{m}(s,t))^2\right\},\end{split}
\end{equation}
where $\Ku_{m}(s,t)$ are the kernels $\Ku_{m}(s,t)=\sum_{k=0}^{m} {{R}_k^{\alpha,\beta}(s)_q
{R}_k^{\alpha,\beta}(t)_q}/{d_k^2}$,
and  $d_n^2$ denotes the squared norm of the $n$-th \textit{standard} $q$-Racah polynomials
(see table \ref{tabla1}). Notice that $\forall A,B>0$ and $a\neq b$, {$\Ku_{m}(a,b)>0$}.
Thus, by the Proposition 1 in \cite{RR}
the polynomials $R_n^{\alpha,\beta,A, B}(s)_q$ are well defined for all values $A,B>0$.

\begin{rem}
If $A,B$ are in general complex numbers then, according to Proposition 1 in \cite{RR}, in 
order that there exists a sequence of orthogonal polynomials 
$(R_n^{\alpha,\beta,A,B}(s)_q)$ the condition 
$\kappa_{n-1}^{\alpha,\beta}(0,N)\neq0$, where $\kappa_{n-1}^{\alpha,\beta}(0,N)$ is defined in 
\refe{K(a,b)}, should be hold for all $n\in\NN$, $A,B\in\CC$.
\end{rem}

\subsection*{Representation formulas for ${R}_n^{\alpha,\beta,A,B}(s)_q$}

Let us now obtain some explicit formulas for the $q$-Racah-Krall polynomials.
The first formula follows from \cite[Eq. (29)]{RR} (see also \refe{repform} from above)
\bq\begin{split}
\label{rez1-rac}
R_n^{\alpha,\beta,A, B}(s)_q =R_n^{\alpha,\beta}(s)_q -
AR_n^{\alpha,\beta,A, B}(0)_q\Ku_{n-1}(s, 0)-
BR_n^{\alpha,\beta,A, B}(N)_q\Ku_{n-1}(s, N)
\end{split}
\eq
where
\begin{equation}\label{ker-a}
\K_{n-1}^{\alpha,\beta}(s, 0)=\varkappa_0^{\alpha,\beta}(s, n)
R_{n-1}^{\alpha,\beta}(s)_q+\conju{\varkappa}_0^{\alpha,\beta}(s, n)
\frac{\nabla R_{n-1}^{\alpha,\beta}(s)_q}{\nabla x(s)},
\end{equation}
\begin{equation}\label{rez17}
\begin{split}
 \varkappa_0^{\alpha,\beta}(s,n)  =\dst & \frac{(\delta \gamma q)^{-n+1}(1-\alpha\beta q^n)(1-\delta \gamma q^{s+n})
(\alpha\beta q,\alpha\beta q^{n+1}; q)_{n-1}}{(1-\alpha\beta q)(1-\delta \gamma q^{s+1})
(q,\alpha\beta \gamma^{-1}q,\alpha\delta^{-1}q,\beta q;q)_{n-1}},\\
\conju{\varkappa}_0^{\alpha,\beta}(s,n) =\dst & \frac{q^{2-s}(\delta \gamma)^{-n+3}
(\alpha\beta q, \alpha\beta q^{n+1};q)_{n-1}(q^s-1)(q^s-\delta^{-1})}
{(1-\alpha\beta q)(1-\delta \gamma q^{2s+2})}\\
&\times\frac{(q^s-\beta\gamma^{-1})(q^s-\alpha\delta^{-1}\gamma^{-1})}{(q,\alpha\beta \gamma^{-1}q,\alpha\delta^{-1}q,\beta q;q)_{n-1}},
\end{split}
\end{equation}
and
\begin{equation}\label{ker-b}
\Ku_{n-1}(s, N)=\varkappa_N^{\alpha,\beta}(s, n)R_{n-1}^{\alpha,\beta}(s)_q+
\conju{\varkappa}_N^{\alpha,\beta}(s, n)
\frac{\Delta R_{n-1}^{\alpha,\beta}(s)_q}{\Delta x(s)},
\end{equation}
\begin{equation*}
\begin{split}
\varkappa_N^{\alpha,\beta}(s, n)= & \dst \frac{(\gamma )^{-n+1}(1-\alpha\beta q^{n})(1-\delta q^{s-n+1})(\alpha\beta q,\alpha\beta q^{n+1};q)_{n-1}}
{(1-\alpha\beta)(1-\delta q^s)(\beta^{-1}\delta^{-1}q^{-n},\alpha q,\beta\delta q,q,\alpha\beta \gamma^{-1}q;q)_{n-1}}, \\
\conju{\varkappa}_N^{\alpha,\beta}(s, n)= & \dst\frac{(\gamma )^{-n+1}q^{-s}
(1-\alpha q^{s+1})(1-\beta\delta q^{s+1})(1-\gamma q^{s+1})(1-\delta \gamma q^{s+1})}
{(1-\alpha\beta q)(1-\delta \gamma q^{2s+2})}\\
&\times\frac{(\alpha\beta q, \alpha\beta q^{n+1};q)_{n-1}}{(\beta^{-1}\delta^{-1}q^{-n},\alpha q,\beta\delta q,q,\alpha\beta \gamma^{-1}q; q)_{n-1}}.
\end{split}
\end{equation*}
Inserting \refe{ker-a} and \refe{ker-b} into the formula \refe{rez1-rac} one finds
the 1{\it st} representation formula
\begin{equation}\label{rez63a}
\begin{split}
\dst{R}_n^{\alpha,\beta,A, B}(s)_q  = \dst  R_n^{\alpha,\beta}(s)_q+
\overline{A}(s, n)R_{n-1}^{\alpha,\beta}(s)_q & +\overline{B}(s, n)
\frac{\nabla R_{n-1}^{\alpha,\beta}(s)_q}{\nabla x(s)}\\  &
+\overline{C}(s, n)\frac{\Delta R_{n-1}^{\alpha,\beta}(s)_q}{\Delta x(s)}
\end{split}
\end{equation}
with
\begin{equation*}
\begin{split}
\overline{A}(s, n)=& -A{R}_n^{\alpha,\beta,A, B}(0)_q\varkappa_0^{\alpha,\beta}(s, n)-
B {R}_n^{\alpha,\beta,A, B}(N)_q\varkappa_N^{\alpha,\beta}(s, n),\\[4mm]
\overline{B}(s, n)=& -A{R}_n^{\alpha,\beta,A, B}(0)_q\conju{\varkappa}_0^{\alpha,\beta}(s, n),
\quad
 \overline{C}(s, n)=-B{R}_n^{\alpha,\beta,A, B}(N)_q\conju{\varkappa}_N^{\alpha,\beta}(s, n),
\end{split}
\end{equation*}
where ${R}_n^{\alpha,\beta,A, B}(0)_q$ and ${R}_n^{\alpha,\beta,A, B}(N)_q$ are given
in \refe{u-u2m}.
Notice that $R_n^{\alpha,\beta,A,B}(s)_q$ defined in \refe{rez1-rac} is a polynomial
of degree $n$ in $x(s)$. However, it is not easy to see that ${R}_n^{\alpha,\beta,A, B}(s)_q$ 
is a polynomial of degree $n$ in $x(s)$ by the formula \refe{rez63a} 
since the functions $\overline{A}$, $\overline{B}$ and $\overline{C}$ as well as 
${\nabla R_{n-1}^{\alpha,\beta}(s)_q}/{\nabla x(s)}$ and 
${\Delta R_{n-1}^{\alpha,\beta}(s)_q}/{\Delta x(s)}$ in \refe{rez63a} are not, in general, 
polynomials in $x(s)$. 

The 2{\it nd} representation formula of the $q$-Racah-Krall polynomials 
follows from \cite[section 3, (17)]{RR} and has the following form
\bq \label{repfor-n-rac}
\phi(s){R}_n^{\alpha,\beta,A, B}(s)_q=A(s,n)R_n^{\alpha,\beta}(s)_q+
B(s,n)R_{n-1}^{\alpha,\beta}(s)_q,
\eq
 \begin{equation}\label{ABsn}
\begin{split}
\phi(s) & =(1-\delta\gamma q^{N+1} q^s)(1-\delta \gamma q^{s+1})(q^{-s}-1)(q^{-s}-q^{-N}), \\[4mm]
A(s, n)& =\phi(s) - \dst\frac{1}{d_{n-1}^2 } \Big\{
A {R}_n^{\alpha,\beta,A, B}(0)_q R_{n-1}^{\alpha,\beta}(0)_q (q^{-s}-q^{-N})(1-\delta\gamma q^{N+1} q^s) \\
&  \qquad\quad +
B {R}_n^{\alpha,\beta,A, B}(N)_q R_{n-1}^{\alpha,\beta}(N)_q(1-\delta \gamma q^{s+1})(q^{-s}-1) \Big\}, \\
B(s, n)& = \dst\frac{1}{d_{n-1}^2 } \Big\{
A {R}_n^{\alpha,\beta,A, B}(0)_q R_{n}^{\alpha,\beta}(0)_q (q^{-s}-q^{-N})(1-\delta\gamma q^{N+1} q^s)\\
&  \qquad\quad +
B {R}_n^{\alpha,\beta,A, B}(N)_q R_{n}^{\alpha,\beta}(N)_q(1-\delta \gamma q^{s+1})(q^{-s}-1)  \Big\},  \\
\end{split}
\end{equation}
where ${R}_n^{\alpha,\beta,A, B}(0)_q$ and ${R}_n^{\alpha,\beta,A, B}(N)_q$ are given
in \refe{u-u2m}. Note that the functions $\phi(s)$ and $A(s,n)$ 
are polynomials of degree $2$ in $x(s)$ and $B(s,n)$ is 1st degree
polynomial in $x(s)$. Therefore it is obvious from the right hand side of the formula 
\refe{repfor-n-rac} that $\phi(s)R_n^{\alpha,\beta,A,B}(s)_q$ is a polynomial of degree $n+2$
in $x(s)$.

Another representation formula for the $q$-Racah-Krall polynomials
may be obtained by putting the relation \refe{ThetaXi} into \refe{repfor-n-rac}
\[
\phi(s){R}_n^{\alpha,\beta,A, B}(s)_q=a(s;n)R_n^{\alpha,\beta}(s)_q+b(s;n)R_n^{\alpha,\beta}(s+1)_q,
\]
where $a(s;n)=A(s;n)+B(s;n)\Theta(s;n)$, $b(s;n)=B(s;n)\Xi(s;n)$,
and $A$, $B$ and $\Theta$, $\Xi$ are given by \refe{ABsn} and \refe{ThetaXi}, respectively.

We note that from Proposition 3 in \cite{RR} it follows that the polynomials ${R}_n^{\alpha,\beta,A, B}(s)_q$ 
satisfy a second order linear difference equation where the coefficients can be computed explicitly.
Since the expression is large enough we will omit them here.

Finally, by use of \cite[Eq (33), section 3]{RR} one can obtain the TTRR for $R_n^{\alpha,\beta,A,B}(s)_q$
$$
x(s)R_n^{\alpha,\beta,A,B}(s)_q= R_{n+1}^{\alpha,\beta,A,B}(s)_q+
\widetilde{\beta}_n  R_{n}^{\alpha,\beta,A,B}(s)_q+
\widetilde{\gamma}_n  R_{n-1}^{\alpha,\beta,A,B}(s)_q
$$
with the following coefficients
\begin{small}
\begin{equation}\label{TTRR}
\begin{split}
\widetilde{\beta}_n &=   \beta_n
-A \left(\frac{{R}_n^{\alpha,\beta,A,B}(0)_q R_{n-1}^{\alpha,\beta}(0)_q}{{d}_{n-1}^2}-
\frac{{R}_{n+1}^{\alpha,\beta,A,B}(0)_qR_{n}^{\alpha,\beta}(0)_q }{{d}_{n}^2}\right) \\
&-B \left(\frac{{R}_n^{\alpha,\beta,A,B}(N)_q R_{n-1}^{\alpha,\beta}(N)_q }{{d}_{n-1}^2}-
\frac{{R}_{n+1}^{\alpha,\beta,A,B}(N)_q R_{n}^{\alpha,\beta}(N)_q }{{d}_{n}^2}\right), \\
\widetilde{\gamma}_n &= \gamma_n\frac{1\!+\!\Delta_n^{\alpha,\beta,A,B}}{1\!+\!\Delta_{n-1}^{\alpha,\beta,A,B}},\,\, \Delta_n^{\alpha,\beta,A,B}=
\frac{ A {R}_n^{\alpha,\beta,A,B}(0)_q R_{n}^{\alpha,\beta}(0)_q }{{d}_{n}^2}\!+\!
\frac{B {R}_n^{\alpha,\beta,A,B}(N)_q R_{n}^{\alpha,\beta}(N)_q}{{d}_{n}^2}
\end{split}
\end{equation}
\end{small}
where we use the notations defined in Eqs. \refe{u-u2m}.

 \subsection*{Representation of ${R}_n^{\alpha,\beta,A,B}(s)_q$ in terms of basic series}

Let us  now introduce the representation of 
${R}_n^{\alpha,\beta,A,B}(s)_q$ in terms of the basic hypergeometric series. To this end,
we substitute \refe{pol-rac} into the formula \refe{repfor-n-rac} which leads to
\[
\ba{rl}
\phi(s){R}_n^{\alpha,\beta,A, B}(s)_q = &
\dst\frac{(\alpha q; q)_{n-1}(\beta\delta q; q)_{n-1}(\gamma q; q)_{n-1}}
{(\alpha\beta q^n; q)_{n-1}}\times\\
& \dst \sum_{k=0}^{\infty}\frac{(q^{-n}, \alpha\beta q^n, q^{-s}, \delta \gamma q^{s+1}; q)_k}{(\alpha q,
\beta\delta q, \gamma q, q; q)_k}q^k
\Pi_1(q^k),
\ea
\]
where $\phi(s)$, $A(s, n)$ and $B(s, n)$ are given in \refe{ABsn} and
\bq\label{Pi_1}
\begin{split}
\Pi_1(q^k)&  = 
A(s, n)\frac{(1-\alpha q^n)(1-\beta\delta q^n)(1-\gamma q^{n})(1-\alpha\beta q^{n+k})}
{(1-\alpha\beta q^{2n-1})(1-\alpha\beta q^{2n})}
+ B(s, n)\frac{(1-q^{-n+k})}{(1-q^{-n})}\\
& = -\frac{1}{1-q^{-n}}\Big\{\alpha\beta q^nA(s, n)
\vartheta^{\alpha,\beta,\delta,\gamma}_n+B(s, n)q^{-n}\Big\}(q^k-q^{\beta_1}),
\end{split}
\eq
where
\[\begin{split}
q^{\beta_1} & =\frac{A(s, n) \vartheta^{\alpha,\beta,\delta,\gamma}_n+ B(s, n)}
{A(s, n)\alpha\beta q^n\vartheta^{\alpha,\beta,\delta,\gamma}_n+ B(s, n)q^{-n}}, \\
\vartheta^{\alpha,\beta,\delta,\gamma}_n & =
\frac{(1-\alpha q^n)(1-\beta\delta q^n)(1-\gamma q^{n})(1-q^{-n})}
{(1-\alpha\beta q^{2n-1})(1-\alpha\beta q^{2n})}.
\end{split}
\]
Moreover, by use of the identity $(q^k-q^{z})(q^{-z};q)_k=
(1-q^{z})(q^{1-z}; q)_k$ we arrive at the following representation
\bq\label{hypreprac2}
\begin{split}\dst
\phi(s){R}_n^{\alpha,\beta,A, B}(s)_q =& {D}_n^{\alpha,\beta,\beta_1}(s)
{}_{5}\varphi_4 \! \left(\!\ba{c} q^{-n},\alpha\beta q^n, q^{-s},
\delta \gamma q^{s+1}, q^{1-\beta_1} \\ \alpha q,\beta\delta q,\gamma q, q^{-\beta_1}\ea\!\!
\bigg|q, q \!\right),
\end{split}
\eq
where
\begin{equation*}
\begin{array}{rl}
{D}_n^{\alpha,\beta,\beta_1}(s)&=\dst-\frac{(\alpha q; q)_{n-1}(\beta\delta q; q)_{n-1}(\gamma q; q)_{n-1}}
{(\alpha\beta q^n; q)_{n-1}}(1-q^{\beta_1})\\
&\dst\times\frac{1}{1-q^{-n}}\Big\{A(s, n) \alpha\beta q^n\vartheta^{\alpha,\beta,\delta,\gamma}_n+
B(s, n)q^{-n}\Big\}.
\end{array}
\end{equation*}

\begin{rem} Notice that $\phi(s){R}_n^{\alpha,\beta,A, B}(s)_q$ in \refe{hypreprac2}
is a polynomial of degree $n+2$ in $x(s)$. 
To see that formula \refe{hypreprac2} gives a polynomial of degree $n+2$ it is 
sufficient to notice that the function $\Pi_1$ defined in \refe{Pi_1} is a polynomial 
of degree $2$ in $x(s)$
since $A(s,n)$ and $B(s,n)$ defined in \refe{ABsn} are {polynomials} of degree 
2 and 1 in $x(s)$, respectively. 
\end{rem}   

Finally let us mention that the direct substitution of \refe{pol-rac} into 
\refe{repfor-n-rac} leads to the following representation formula
\[
\begin{split}
 \phi(s){R}_n^{\alpha,\beta,A,B}(s)_q&= A(s,n){\Lambda}_n^{\alpha,\beta,\delta,\gamma}
{}_{4}\varphi_3 \left(\ba{c} q^{-n},\alpha\beta q^{n+1}, q^{-s},
\delta\gamma q^{s+1}  \\ \alpha q,\beta\delta q,\gamma q  \ea
\,\bigg|\, q \,,\, q \right)\\
&+B(s, n){\Lambda}_{n-1}^{\alpha,\beta,\delta,\gamma}
{}_{4}\varphi_3 \left(\ba{c} q^{-n+1},\alpha\beta q^{n}, q^{-s},
\delta \gamma q^{s+1}  \\ \alpha q,\beta\delta q,\gamma q  \ea
\,\bigg|\, q \,,\, q \right),
\end{split}
\]
where
$$
{\Lambda}_n^{\alpha,\beta,\delta,\gamma}=\frac{(\alpha q,\beta\delta q, \gamma q;q)_n}
 {(\alpha\beta q^n;q)_n}.
$$

\subsection{The case of one mass point}
In this section we introduce the \textit{standard} $q$-Racah polynomials but with one mass point at the value
$s=0$. All the formulas follow by replacing $B=0$ into the ones in section \ref{2masspoint}.
For example the first representation formula of ${R}_n^{\alpha,\beta,A}(s)_q $ is produced 
by inserting $B=0$ into \refe{rez63a}, 
\begin{equation*}
\begin{array}{rl}\dst
\dst{R}_n^{\alpha,\beta,A}(s)_q & = \dst  R_n^{\alpha,\beta}(s)_q+
\overline{A}(s, n)R_{n-1}^{\alpha,\beta}(s)_q +\overline{B}(s, n)
\frac{\nabla R_{n-1}^{\alpha,\beta}(s)_q}{\nabla x(s)},
\end{array}
\end{equation*}
where 
\begin{equation*}
\begin{array}{rl}
&\overline{A}(s, n)=-A{R}_n^{\alpha,\beta,A}(0)_q\varkappa_0^{\alpha,\beta}(s, n),
\quad\overline{B}(s, n)=-A{R}_n^{\alpha,\beta,A}(0)_q
\conju{\varkappa}_0^{\alpha,\beta}(s, n),
\end{array}
\end{equation*}
and
\begin{equation}\label{u-u1m}
{R}_n^{\alpha,\beta,A}(0)_q=\frac{R_n^{\alpha,\beta}(0)_q}{1+A\Ku_{n-1}(0, 0)}
\end{equation}
in which $\varkappa_0^{\alpha,\beta}(s, n)$ and $\conju{\varkappa}_0^{\alpha,\beta}(s, n)$
are defined in \refe{rez17},
respectively. Moreover, the second representation formula follows
by evaluating \refe{repfor-n-rac} for $B=0$
\[
\phi(s){R}_n^{\alpha,\beta,A}(s)_q=A(s;n)R_n^{\alpha,\beta}(s)_q+
B(s;n)R_{n-1}^{\alpha,\beta}(s)_q,
\]
where $\phi(s)=(q^{-s}-1)(1-\delta\gamma q^{s+1})$,
\begin{equation}\label{ABsn-1m}
\begin{split}
\quad A(s, n) & =\phi(s) - \dst\frac{A}{d_{n-1}^2 } {R}_n^{\alpha,\beta,A}(0)_q R_{n-1}^{\alpha,\beta}(0)_q , \quad
B(s, n)  = \dst\frac{A}{d_{n-1}^2 } {R}_n^{\alpha,\beta,A}(0)_q R_{n}^{\alpha,\beta}(0)_q,  \\
\end{split}
\end{equation}
and ${R}_n^{\alpha,\beta,A}(0)_q$ is  given in \refe{u-u1m}. Finally, 
the third representation formula follows by use of the same idea
\bq\label{repfor-n-rac-1m}
\phi(s){R}_n^{\alpha,\beta,A}(s)_q=a(s;n)R_n^{\alpha,\beta}(s)_q+b(s;n)R_n^{\alpha,\beta}(s+1)_q,
\eq
where $a(s;n)=A(s;n)+B(s;n)\Theta(s;n)$ and $b(s;n)=B(s;n)\Xi(s;n)$, 
and $A$, $B$ and $\Theta$, $\Xi$ are given by \refe{ABsn-1m} and \refe{ThetaXi}, respectively.
Notice that, as for two mass points case, from the above representation formula \refe{repfor-n-rac-1m}
the second order difference equation for ${R}_n^{\alpha,\beta,A}(s)_q$ follows \cite{R,RR}.

Furthermore, one can obtain the coefficients of the TTRR 
by replacing $B=0$ into \refe{TTRR} as the following
 \begin{equation*}
\begin{split}
\widetilde{\beta}_n  & = \beta_n
-A \left(\frac{{R}_n^{\alpha,\beta,A}(0)_q R_{n-1}^{\alpha,\beta}(0)_q}{{d}_{n-1}^2}-
\frac{{R}_{n+1}^{\alpha,\beta,A}(0)_qR_{n}^{\alpha,\beta}(0)_q }{{d}_{n}^2}\right), \\
\widetilde{\gamma}_n &= \gamma_n\frac{1+\Delta_n^{\alpha,\beta,A}}{1+\Delta_{n-1}^{\alpha,\beta,A}},\quad \Delta_n^{\alpha,\beta,A}=
\frac{ A {R}_n^{\alpha,\beta,A}(0)_q R_{n}^{\alpha,\beta}(0)_q }{{d}_{n}^2}.
\end{split}
\end{equation*}

Notice that putting $B=0$ in the basic series representation formulas
\refe{hypreprac2} we obtain the corresponding
basic series representations for the $q$-Racah-Krall polynomials ${R}_n^{\alpha,\beta,A}(s)_q$.


\subsection{Some limit cases}
We first consider the modification of \textit{standard} dual $q$-Hahn polynomials defined 
on the lattice $x(s)=q^{-s}+\gamma\delta q^{s+1}$ by
\[
\begin{split}\dst
R_n(x(s),\gamma, \delta, N|q) & = \dst
(q^{-N}, \gamma q; q)_n 
{}_{3}\varphi_2 \left(\ba{c} q^{-n},  q^{-s},
\delta \gamma q^{s+1}  \\ q^{-N}, \gamma q  \ea
\bigg| q ; q \right),
\end{split}
\]
that are related with the $q$-Racah polynomials by the expression  \cite{ks}
\bq\label{rac->dual}
\lim_{\beta\to 0} R_n^{\alpha,\beta}(s)_q={R}_n(s)_q.
\eq

In order to obtain the modification of \textit{standard} dual $q$-Hahn polynomials
by adding two mass points at the end of the interval of orthogonality, i.e.,
$R_n^{A,B}(s)_q$ satisfying the orthogonality relation
\bq\label{ort-dual}
\begin{split}
\sum_{s = 0 }^{N} R_n^{A,B}(s)_q R_m^{A,B}(s)_q
\rho(s)& \Delta x(s-1/2) + A R_n^{A,B}(0)_q
R_m^{A,B}(0)_q\\& 
+BR_n^{A,B}(N)_qR_m^{A,B}(N)_q = \delta_{n,m} \widetilde{d}_n^2.
\end{split}
\eq
where 
$$\rho(s)=
\frac{q^{Ns-(^s_2)}(\gamma q)^N}{(-\gamma)^s(1-\gamma\delta q)(q^{-1/2}-q^{1/2})}
\frac{(\delta q; q)_{N}}{(\gamma\delta q^2; q)_N}
\frac{(\gamma q,\gamma\delta q,q^{-N}; q)_{s}}
{(q, \delta \gamma q^{N+2}, \delta q; q)_s}.
$$
Using the same procedure as before it follows that
\begin{small}
\begin{equation*}
\begin{array}{l}\dst
R_n^{A,B}(0)_q=\frac{(1+B \K_{n-1}(N, N))R_n(0)_q-
B\K_{n-1}(0, N)R_n(N)_q}{\kappa_{n-1}(0,N)},\\
R_n^{A,B}(N)_q=\dst
\frac{-A\K_{n-1}(N, 0)R_n(0)_q+(1+A\K_{n-1}(0, 0))R_n(N)_q}{\kappa_{n-1}(0,N)},
\end{array}
\end{equation*}
\end{small}%
\begin{small}
\begin{equation*}
\begin{split}
\widetilde{d}_n^2  =
d_n^2 & \!+\! \frac{A (R_n(0)_q)^2 \{1+ B\K_{n-1}(N,N)\} \! +\!
B(R_n(N)_q)^2 \{ 1\!+\! A\K_{n-1}(0,0) \} }{\kappa_{n-1}(0, N)} \\
& -\frac{2 AB R_n(0)_q R_n(N)_q\K_{n-1}(0,N)}
{\kappa_{n-1}(0, N)}.
\end{split}
\end{equation*}
\end{small}%
In the above formulas
\begin{equation*}
\kappa_{m}(s,t)=1+A\K_{m}(s,s)+ B\K_{m}(t,t)
+AB\left\{\K_{m}(s,s)\K_{m}(t,t)-(\K_{m}(s,t))^2\right\},
\end{equation*}
$\K_m(s,t)=\sum_{k=0}^{m} {R_n(s)_q R_n(t)_q}/{d_k^2}$,
and  $d_n^2$ denote $m$-th Kernel and the norm of the \textit{standard} dual $q$-Hahn polynomials
$$
d_n^2=(\gamma\delta q)^n(q,\delta^{-1}q^{-N},\gamma q,q^{-N};q)_n.
$$

If we now fix $\alpha q=q^{-N}$ in the orthogonality relation
for the \textit{standard} $q$-Racah polynomials \refe{ort-umod} and take the limit
$\beta\to 0$, then using \refe{rac->dual} we obtain the orthogonality relation \refe{ort-dual}, 
and therefore, it is straightforward to see that
\[
\lim_{\beta\to 0} R_n^{\alpha,\beta,A,B}(s)_q={R}_n^{A,B}(s)_q.
\]
Thus, all properties of the modified dual $q$-Hahn polynomials
 ${R}_n^{A, B}(s)_q$ can be obtained from the corresponding
properties of the modified $q$-Racah polynomials
$R_n^{\alpha,\beta,A,B}(s)_q$ by taking the appropriate limit. 

To conclude this paper we discuss two other important limit cases of the
$q$-Racah-Krall polynomials ${R}_n^{\alpha, \beta, A,B}(s)_q$.

The first one is when we take the limit $q\to1$. In fact, as
$q\to1$ in \refe{pol-rac} we recover the Racah polynomials in the quadratic lattice 
$x(s)=s(s+\gamma+\delta+1)$ \cite{ks}, i.e.,
$$
\lim_{q\to1} {R}_n(q^{-s}+\delta\gamma q^{s+1},q^{\alpha},q^{\beta},q^{\gamma},q^{\delta})_q= 
{R}_n(s(s+\gamma+\delta+1),\alpha, \beta,\gamma,\delta),
$$
where
\[
\begin{split}\dst
R_n^{\alpha, \beta}(s)= & \dst
\frac{(\alpha+1)_n (\beta+\delta+1)_n
 (\gamma+1)_n}{(\alpha+\beta+n+1)_n}\times \\ & 
{}_{4}\mbox{F}_3 \left(\ba{c} {-n}, { \alpha+\beta+n+1}, {-s},
{s+\gamma+\delta+1}  \\ {\alpha+1}, {\beta+\delta+1}, { \gamma+1}  \ea
\,\bigg| \,1\, \right).
\end{split}
\]
Straightforward calculations  show that all the properties of the \textit{standard} $q$-Racah 
polynomials ${R}_n^{\alpha, \beta}(s)_q$ becomes into the properties of the 
\textit{standard} Racah ones. In particular, we have the following limit relation
$$
\lim_{q\to1} {R}_n^{\alpha, \beta,A,B}(q^{-s}+\delta\gamma q^{s+1})_q= {R}_n^{\alpha, \beta,A,B}(s(s+\delta+\gamma+1)),
$$
where ${R}_n^{\alpha, \beta,A,B}(s(s+\delta+\gamma+1))$ denotes the 
modification of the Racah polynomials by adding 
two delta Dirac masses at the points $0$ and $N$.
Thus, taking appropriate limits one can construct the analogue of the
Askey Tableau but for the Krall type polynomials.

Another important family of Krall-type polynomials are the so called $q$-Hahn-Krall tableau of
orthogonal polynomials considered in \cite{RJ, RR}. In order to obtain it
first of all notice that the \textit{standard} $q$-Racah  polynomials 
are defined on the lattice $x(s)=c_1q^{s}+c_2q^{-s}+c_3$
where $c_1=\gamma\delta q$, $c_2=1$ and $c_3=0$. 
Then, making the transformation  {$q^{\delta-N}x(s)=c_1\mu(s)$}, $\alpha=\nu$,
$\beta=\mu$, $\gamma q\to q^{-N}$ and $\delta\to q^{\delta}$ in \refe{pol-rac} and then taking
the limit $q^{\delta}\rightarrow 0$ by use of the identity \cite{ran}
\begin{small}
$$
(q^{-s}; q)_k(q^{s+\zeta}; q)_k\!=\!(-1)^kq^{k(\zeta+\frac{k-1}{2})}\prod_{i=0}^{k-1}
\Big[\frac{x(s)\!-\!c_3}{c_1}-q^{-\frac{\zeta}{2}}(q^{i+\frac{\zeta}{2}}+q^{-i-\frac{\zeta}{2}})\Big],
$$
\end{small}%
we obtain
\bq
\dst
C_nR_n(x(s),\alpha,\beta,\gamma,\delta|q) \stackrel{q^{\delta}\rightarrow 0}{\longrightarrow}  h_n^{\mu, \nu}(\mu(s); N|q),
\label{limit-rac-han-nu}
\eq
where $h_n^{\mu, \nu}(\mu(s); N|q)$ are the $q$-Hahn polynomials on the lattice {$\mu(s)=q^{-s}$ }
$$
\dst h_n^{\mu, \nu}(\mu(s); N|q):=
\frac{(\nu q;q)_n (q^{-N};q)_n}{(\mu\nu q^{n+1};q)_n}
{}_{3}\varphi_2 \left(\ba{c} q^{-n},\mu\nu q^{n+1}, \mu(s)  \\ \nu q,q^{-N} \ea
\,\bigg|\, q \,,\, q \right),
$$
and $C_n:=1$. 
Then we get the $q$-Hahn-Krall polynomials by use of the above limit relation and
the property that
\bq\label{limit-rac-han-pts}
C_nR_n^{\alpha, \beta}(0)\stackrel{q^{\delta}\rightarrow 0}{\longrightarrow} h_n^{\mu, \nu}(\mu(0); N|q)\quad
C_nR_n^{\alpha, \beta}(N)\stackrel{q^{\delta}\rightarrow 0}{\longrightarrow} h_n^{\mu, \nu}(\mu(N); N|q),
\eq
and
$$
C_kd_k^2\stackrel{q^{\delta}\rightarrow 0}{\longrightarrow} \overline{d}_k^2=
(-\nu q)^nq^{(^n_2)-Nn}\frac{(q, \mu q, \nu q, q^{-N}, \mu\nu q^{N+2}; q)_n}{(\mu\nu q^2; q)_{2n}
(\mu\nu q^{n+1}; q)_n},
$$
where $d_k$ and  $\overline{d}_k ^2$ denote the norms for the \textit{standard} $q$-Racah and 
the $q$-Hahn polynomials, respectively.
The following limit relation for the kernels of the \textit{standard} $q$-Racah and $q$-Hahn polynomials
are introduced by applying the aforesaid transformation 
\bq\label{limit-rac-han-kernel}\begin{split}
& \K_n^{\alpha, \beta}(s_1, s_2):= 
\dst\sum_{k=0}^{n}{\dst\frac{C_kR_k^{\alpha, \beta}(s_1)_qC_kR_k^{\alpha, \beta}(s_2)_q}{C_k^2d_k^2}}
\stackrel{q^{\delta}\rightarrow 0}{\longrightarrow} \\
& \qquad \dst\sum_{k=0}^{n}{\dst\frac{h_k^{\mu, \nu}(x(\overline{s}_1); N|q)
h_k^{\mu, \nu}(x(\overline{s}_2); N|q)}{\overline{d}_k^2}}:=
\K_n^{\mu, \nu}(\overline{s}_1, \overline{s}_2).
\end{split}
\eq
Therefore, as a result of the limit relations defined by \refe{limit-rac-han-nu}, \refe{limit-rac-han-pts}, and \refe{limit-rac-han-kernel} and the formula by \refe{rez1-rac}, we obtain that
$$
\lim_{q^a\rightarrow 0}
C_n R_n^{\alpha, \beta, A, B}(s)_q =  h_n^{\mu, \nu, A, B}(\mu(s);N|q):= h_n^{\mu, \nu, A, B}(s)_q,
$$
where $C_n=1$. In other words, we obtain the $q$-Hahn-Krall polynomials on the lattice {$\mu(s)=q^{-s}$}
which satisfy the orthogonality relation
\[\begin{split}
\sum_{s = 0 }^{N}h_n^{\mu, \nu, A, B}(s)_q & h_m^{\mu, \nu, A, B}(s)_q\rho(s)
\Delta x(s-\mbox{$\frac12$})
+Ah_n^{\mu, \nu, A, B}(0)_qh_m^{\mu, \nu, A, B}(0)_q\\
&+Bh_n^{\mu, \nu, A, B}(N)_q
h_m^{\mu, \nu, A, B}(N)_q=\delta_{n,m}\overline{d}_k^2,\quad {\mu(s)=q^{-s}},
\end{split}
\]
where $\rho$ is the weight function of the $q$-Hahn polynomials \cite[page 445]{ks}.
%
%

\subsection*{Acknowledgements:} This work was partially supported by MTM2009-12740-C03-02
(Ministerio de Econom\'\i a y Competitividad), FQM-262, FQM-4643, FQM-7276 (Junta de Andaluc\'\i a),
Feder Funds (European Union).
The second author is supported by a grant from T\"{U}B\.{I}TAK, the Scientific and
Technological Research Council of Turkey. This research has been done during the stay
of the second author at the Universidad de Sevilla. She also
thanks to the Departamento de An\'alisis Matem\'atico of the Universidad de Sevilla
and IMUS for their kind hospitality. 



\end{document}